\documentclass[review]{elsarticle}
\usepackage{amsfonts,amssymb,amsmath,amsbsy,amsthm,amscd}
\usepackage{lineno,hyperref}
\usepackage[english]{babel}
\usepackage{graphicx}
\usepackage{color}

\newtheorem{theorem}{Theorem}[section]

\newtheorem{lemma}[theorem]{Lemma}

\newtheorem{proposition}{Proposition}

\newtheorem{remark}{Remark}

\renewcommand{\le}{\leqslant}
\renewcommand{\ge}{\geqslant}

\begin{document}

\begin{frontmatter}

\title{The repulsive Euler-Poisson equations with  variable doping profile
}


\author{Olga S. Rozanova}
\ead{rozanova@mech.math.msu.su} 

\address[1]{Department of Mechanics and Mathematrics, Moscow
State University, Moscow 119991 Russia}


\begin{abstract}
 We prove that arbitrary smooth perturbations of the zero equilibrium state of the repulsive pressureless Euler-Poisson  equations, which describe the behavior of cold plasma, blow up for any non-constant doping profile already in one-dimensional space. Further, we study small perturbations of the equilibrium to determine which  properties of the doping profile contribute to the blow-up. We also propose a numerical procedure that allows one to find the blow-up time for any initial data and present examples of such calculations for various doping profiles for  standard initial data, corresponding to the laser pulse.
\end{abstract}


\def\sign{\mathop{\rm sgn}\nolimits}





\begin{keyword}
Euler-Poisson system \sep
equations of cold plasma  \sep singularity formation \sep variable doping profile

\MSC 35F55 \sep 35Q60 \sep 35B44 \sep 35B20  \sep    35L80
\end{keyword}

\end{frontmatter}

\section{Introduction}

The system of Euler-Poisson equations describing the behavior of cold plasma in $ {\mathbb R}^{ d}$ in the repulsive case has the following form:
\begin{equation}\label{1}
V_t+ (V \cdot \nabla ) V =-\nabla \Psi,\quad \dfrac{\partial n}{\partial t}+ {\rm div}\,( n V)
=0,\quad \Delta \Psi =c(x) -n.
\end{equation}
The components of the solution $V$, $n> 0$, $\Psi$ are the velocity, electron density and electric field potential, respectively, they depend on the time $t\ge 0$ and the point $x\in {\mathbb R}^{ d} $. A fixed $C^1$ - smooth function
$c(x)>0$ is the density background or the so-called doping profile.

An extensive literature is devoted to the Euler-Poisson equations. They consider cases of repulsive and attractive forces, zero and non-zero background density. The behavior of the solution in each of these cases is significantly different. The case of repulsive force arises when modeling the phenomena of plasma physics and semiconductors, that is, a medium consisting of electrons. The attractive case refers to a medium in which the gravitational force acts between particles, that is, to astrophysical models (system \eqref{1} corresponds to this case after changing the sign in front of $\Psi$ in the first equation). The Euler-Poisson system was considered both with and without pressure. In all cases, there exist results  about the local existence of a solution to the Cauchy problem; it is known that  there are initial data at which the solution to the Cauchy problem loses its original smoothness \cite{Yuen}; the possibility of constructing a global weak solution was investigated (we refer to  \cite{Chen23} for a recent review). An interesting topic is the study of critical thresholds, that is, the most accurate separation of the initial Cauchy data into two classes, one of which includes those that correspond to globally smooth solutions, and the other to solutions that develop a singularity over a finite time. This work began in \cite{ELT}, where criteria were derived for basic one-dimensional and some more complex situations involving the pressureless Euler-Poisson equations. The current state of the problem can be found in \cite{Bhat23}.

In this paper we consider repulsive Euler-Poisson equations  with a non-zero constant density background describing cold plasma \cite{GR75}. Recently, interest in such models has increased significantly due to the possibility of creating  accelerators on the wake wave \cite{esarey09}. This is the most difficult case of the Euler-Poisson equations, since the solutions are oscillating and even the zero equilibrium  is usually unstable due to the phenomenon of nonlinear resonance.

Since cold plasma is a very unstable medium, from a practical point of view, questions about the possibility of choosing a regime that would guarantee a smooth solution for as long as possible are relevant. However, in the space of many spatial variables ${d}\ge 2$, ${d}\ne 4$, all physically reasonable solutions lose smoothness in a finite time \cite{R22_Rad}, \cite{RT}.
However, for $d=1$ and a constant background density $c(x)=C$, there is a sufficiently large neighborhood of the zero equilibrium $V=V_x=\Phi_x=0$, $n=C$,  such that the solution with data from this neighborhood remains smooth for all $t>0$ \cite{ELT}, \cite{RChZAMP21} (in the case of pressure this also occurs \cite{Guo17}). It was recently shown that for the exceptional dimension $d=4$ there is also a neighborhood of the zero equilibrium corresponding to a globally smooth solution.

 In this paper, we show that this phenomenon does not persist under any non-constant background, including an arbitrarily small perturbation of a constant one. In other words, if there is a point $x_0\in \mathbb R$ such that $c'(x_0)\ne 0$, then any solution whose data are an arbitrarily small perturbation of the equilibrium $V=0$, $ \Psi = \rm const$ blows up in a finite time.

The paper is organized as follows. In Sec.\ref{S2} we rewrite \eqref{1} in terms of velocity and electric field, which allows us to consider the dynamics of solutions along characteristics and obtain a system of ODEs for this purpose.
In Sec.\ref{S3} we prove a rather general lemma that allows us to show that any nontrivial solution to this system blows up in finite time. However, this result does not answer the question about the qualitative characteristics that promote the blow-up. To study this issue, in Sec.\ref{S4} we  linearize the above system of ODEs using Radon's lemma and reduce the blow-up problem to the problem of vanishing of a special scalar function, which obeys a third order linear equation with time-dependent coefficients. This problem can be studied numerically, and in Sec.\ref{S6} we give examples of computation of such kind for several doping profiles.
In Sec.\ref{S5} we consider a perturbation of the zero equilibrium and choose the size of the perturbation as a small parameter. First we study the change of period of oscillations depending on the properties of the doping profile.  Then we use the Floquet theory  to provide independent proof that the solution of  the Cauchy problem with arbitrary small data blow up in a finite time due to nonlinear resonance.
 Sec.\ref{S7} is devoted to discussion. In particular, we outline a method for proving that, in the case of a damped Euler-Poisson system, there exists a neighborhood of the zero equilibrium such that if the initial data together with its derivatives are sufficiently small, then the corresponding solution is globally smooth in time and tends to the zero equilibrium at $ t\to\infty$.

\section{Analysis along characteristics}\label{S2}

We introduce the electric field vector $E=\nabla \Psi$ and
   rewrite \eqref{1} in terms of $V$ and $E$, see, e.g. \cite{CH18} for details.  In the case ${d}=1$ this results in
\begin{equation}\label{2}
\dfrac{\partial V}{\partial t}+ V  \dfrac{\partial V}{\partial x}=-E, \quad \dfrac{\partial E}{\partial t}+ V  \dfrac{\partial
E}{\partial x}=c(x) V, \quad n=c(x)- \dfrac{\partial E}{\partial x}.
\end{equation}
We consider \eqref{2}  with
the Cauchy data
\begin{equation}\label{CD}
(V,E)|_{t=0}=(V_0(x), E_0(x))\in {C^2} ({\mathbb R}).
\end{equation}
The problem \eqref{2}, \eqref{CD} has a solution locally in time that is as smooth as the initial data, and the formation of a singularity is associated with the blow-up of the ether of the solution components themselves or their first derivatives, e.g.
\cite{RY}, \cite{Bhat20}.

Thus, along characteristics $x=x(t)$, starting from a point $x_0\in \mathbb R$ the solution $(V(x(t)),E(x(t)))$ obeys the system of ODEs
\begin{equation}\label{char_sol}
\dot x=V, \quad \dot V =-E, \quad \dot E  =c(x) V,
\end{equation}
with the initial data
\begin{equation}\label{char_sol_CD}
x(0)=x_0, \quad V(0)=V_0(x_0), \quad E(0)=E_0(x_0).
\end{equation}

\section{Nonexistence result for non-constant doping profile}\label{S3}

\begin{theorem}\label{T1}
If $c(x)\ne \rm const$, then any nontrivial solution to the Cauchy problem \eqref{2}, \eqref{CD} blows up in a finite time.
\end{theorem}

First of all, we  prove a lemma that can be applied in many similar situations. In fact, this is a more general formulation of Lemma 2.2. \cite{Carillo}.

\begin{lemma}\label{Lem1}
Assume that the mapping $x\mapsto X(t)$ $(\mathbb R \mapsto \mathbb R)$ is continuous, and the trajectory $X(t)$ is (nontrivially) periodic with respect to $t\ge 0$ for all $x_0\in \mathbb R$ with  period $T(x_0)$ that  depends continuously on $x_0$. Then, if $T(x)$ is not constant, there exist $x_1$ and $x_2$ from $\mathbb R$ such that  $X_1(t_*)=X_2(t_*)$ for some $t_*>0$.
\end{lemma}

 {\it Proof of Lemma \ref{Lem1}} closely follows the proof of density of the orbit of irrational rotation of a circle. Let us fix $x_1$.
Since $X_i(x)$ is periodic, there exist $x_i^-$ and $x_i^+$, such that $x_i^-\le x_i \le x_i^+$ and $X_i(t)\in {\bf X}_i=[x_i^-, x_i^+] $, $i=1,2$. Due to the continuous dependence on the initial point, for sufficiently small $|x_1-x_2|$ the segments ${\bf X}_1$ and ${\bf X}_2$ intersect,  $x_1, x_2 \in {\bf X}={\bf X}_1\cap {\bf X}_2$.

Assume that $T_1=T(x_1)\ne T_2=T(x_2)$, and $\frac{T_1}{T_2}\ne \mathbb Q$ (due to the continuous dependence $T(x)$ on $x $ we can always choose $x_2$ with this property). Let us consider the trajectory $X_2(t)$ on a circle $S_{T_2}$ of length $T_2$. If we prove that the orbit of point $x_2$ constructed as $f_n=X_2(n T_1)$, $n\in \mathbb N$, is dense in ${\bf X}_2$ (and therefore in $\bf X$), we prove by continuity that there exists $t_*>0$ such that $X_1(t_*)=X_2(t_*)$. Indeed, $X_1(n T_1)=x_1$, and the points of the orbit $f_n$ can be found in an arbitrary small neighborhood of $x_1\in {\bf X}_1$.

Let us prove first that the orbit $P_n=x_2+ n T_1\, ({\rm mod} \, T_2 )$ of point $x_2$ is dense in  $S_{T_2}$. We have to show that $P_n$ comes into any arc of length $\epsilon>0$ of the  circle $S_{T_2}$. Let $N\in \mathbb N$ be sufficiently large such that $\frac{1}{N}<\epsilon$. Without loss of generality we can take $x_2=0$. We consider a set of points  $ (0, a, 2a, \dots, Na) ({\rm mod}\, T_2)$, $a=\frac{T_1}{T_2}$. At least two of these $N+1$ points lay on the same arc of length
$\frac{1}{N}$, therefore for some $k$ and $l $ we have
$|k a- la| ({\rm mod}\, T_2) <\epsilon$. Thus, $P_{k-l}$ is the turn of the circle $S_{T_2}$ to the arc $|(k - l)a| ({\rm mod}\, T_2)<\epsilon$ and the points of the sequence
 $(0, (k - l)a, 2(k -l)a, \dots  ) ({\rm mod}\, T_2)$ visit any arc of length $\epsilon$.

 The continuity of $X_2(t)$ implies that $f_n$ is dense in ${\bf X}_2$.
$\Box$

\medskip

 {\it Proof of Theorem \ref{T1}. }

 Let us obtain the equation  for characteristics defined by the system \eqref{char_sol}. Denote
 \begin{equation*}
  C(x)=\int\limits_{x_0}^x \, c(\xi)\, d\xi> 0.
 \end{equation*}
 Then from \eqref{char_sol} we have
 \begin{equation*}
 \dot E=c(x) V = c(x)\dot x=\dot C(x), \quad \dddot x = - \dot E =\dot C(x).
 \end{equation*}
  The latter equation implies
 \begin{equation}\label{xC}
   \ddot x +\mathcal C(x)=0,\quad \mathcal C(x)=C(x)-C(x_0)+E(x_0),
 \end{equation}
 point $x_0$ is the equilibrium, the center. As follows from Lemma \ref{Lem1}, if $x_0$ is not an isochronous center \eqref{xC} (the oscillation period depends on $x_0$), then the characteristics of the system \eqref{char_sol} necessarily intersect.

We use the following theorem \cite{Sabatini} about the properties of
a second order differential equation of Lienard type
\begin{equation}\label{Lienard}
\ddot y+ f (y) \dot y+g(y)=0.
\end{equation}
\begin{theorem}\label{S}
Let $f, g$  be analytic, $g$ odd, $f (0)=g(0)=0$, $g'(0)>0$. Then
$\mathcal O =(y,\dot y)=(0,0)$ is a center if and only if $f$ is odd and
$\mathcal O )$ is an isochronous center for \eqref{Lienard} if and only if
\begin{equation}\label{tau}
\tau(y) :=\left(\int\limits_0^y sf(s) ds \right)^2-y^3 (g(y)-g'(0)y)=0.
\end{equation}
\end{theorem}

To apply this theorem to our case  \eqref{xC}, we change $y=x-x_0$ and move $x_0$ to the origin. Thus, we have
\begin{equation}\label{tauV}
\tau(y) =-y^3 (\mathcal C(y)-c(x_0) y).
\end{equation}
We can see from \eqref{tauV} that
$\tau=0$ if and only if $\mathcal C(y)=c(x_0) y$, therefore $C'(x)=c(x)=c(x_0)=\rm const$.
Thus, a constant $c(x)$ is the only way to satisfy condition \eqref{tau}.
Theorem \ref{T1} is proved. $\Box$



\section{Behavior of derivatives}\label{S4}

Next, let us denote $V_x=v$, $E_x=e$. Then, differentiating \eqref{2} with respect to $x$, we obtain the following system along the characteristic $x=x(t)$:
\begin{equation}\label{char_der}
\dot v=-v^2-e, \quad \dot e =-sv+c(x) v + c'(x) V,
\end{equation}
with the  initial data
\begin{equation}\label{char_der_CD}
 v(0)=V'_0(x_0), \quad s(0)=E'_0(x_0).
\end{equation}

We  see that if $c(x)=C=\rm const$, then system \eqref{char_der} splits from \eqref{char_sol} and can be easily integrated. In this way, it is possible to obtain an analytical criterion for the formation of a singularity in terms of the initial data  \cite{RChZAMP21}, \cite{ELT}.
It implies that if the initial data \eqref{CD} are such that for any $x_0\in \mathbb R$ the inequality
\begin{eqnarray} \label {crit2}
\left (V'_0 (x_0) \right) ^ 2 + 2 \, E'_0 (x_0) -C <0
\end {eqnarray}
holds, then the solution is periodic in $t$ and retains initial smoothness for all $t>0$. We see that for $C>0$ there is a neighborhood of equilibrium $V=E=0$ in the $C^1$ - norm such that a solution starting from this neighborhood preserves smoothness.

Obviously, if $x_0\in \mathbb R$ is such that $c'(x_0)=0$, then from the previous consideration it follows that the solution does not blow up along the specific characteristic $x(t)=V$, $x(0)=x_0$.

In the case $c'(x_0)\ne 0$ we can consider \eqref{char_der} as a system for $v$ and $s$
with time-dependent coefficients known from \eqref{char_sol}.
Using the Radon theorem \cite{Radon}, \cite{Riccati}, we can linearize this system. For convenience, we present the statement of this theorem here.

\begin{theorem}\label{TR} [The Radon lemma]
\label{T2} A matrix Riccati equation
\begin{equation}
\label{Ric}
 \dot W =M_{21}(t) +M_{22}(t)  W - W M_{11}(t) - W M_{12}(t) W,
\end{equation}
 {\rm (}$W=W(t)$ is an $(n\times m)$ matrix , $M_{21}$ is an $(n\times m)$ matrix , $M_{22}$ is an $(m\times m)$ matrix, $M_{11}$ is an
 $(n\times n)$ matrix, $M_{12} $ is an $(m\times n)$ matrix{\rm )} is equivalent to the homogeneous linear matrix equation
\begin{equation}
\label{Lin}
 \dot Y =M(t) Y, \quad M=\left(\begin{array}{cc}M_{11}
 & M_{12}\\ M_{21}
 & M_{22}
  \end{array}\right),
\end{equation}
 {\rm (}$Y=Y(t)$  is an $(n\times (n+m))$ matrix, $M$ is an $((n+m)\times (n+m))$  matrix{\rm )} in the following sense.

Let on some interval ${\mathcal J} \in \mathbb R$ the
matrix-function $\,Y(t)=\left(\begin{array}{c}{Q}(t)\\ {P}(t)
  \end{array}\right)$ {\rm (}${Q}$  is an $(n\times n)$ matrix, ${P}$  is an $(n\times m)$ matrix{\rm ) } be a solution of \eqref{Lin}
  with the initial data
  \begin{equation*}\label{LinID}
  Y(0)=\left(\begin{array}{c}{\mathbb I}\\ W_0
  \end{array}\right)
  \end{equation*}
   {\rm (}$ \mathbb I $ is the identity $(n\times n)$ matrix, $W_0$ is a constant $(n\times m)$ matrix{\rm ) } and  $\det {Q}\ne 0$ on ${\mathcal J}$.
  Then
{\bf $ W(t)={P}(t) {Q}^{-1}(t)$} is the solution of \eqref{Ric}   on ${\mathcal J}$   with
$W(0)=W_0$.
\end{theorem}

\medskip

We see that system  \eqref{char_der} can be written as \eqref{Ric} with
\begin{eqnarray*}\label{M}
W=\begin{pmatrix}
  v\\
  e
\end{pmatrix},\quad
M_{11}=\begin{pmatrix}
  0\\
\end{pmatrix},\quad
 M_{12}=\begin{pmatrix}
  1 & 0\\
\end{pmatrix},\\
M_{21}=\begin{pmatrix}
 0 \\ c'(x(t)) V(x(t))
\end{pmatrix},\quad
M_{22}=\begin{pmatrix}
  0 & -1\\
  c(x(t))  & 0\\
\end{pmatrix}.\\\nonumber
\end{eqnarray*}
Then, according Theorem \ref{TR} the solition of \eqref{char_der} is
\begin{eqnarray*}\
   \quad W(t)=\frac{P(t)}{Q(t)},
 \end{eqnarray*}
 where $P(t)=(p_1(t),p_2(t))^T$ and  $Q(t)$ solves the linear system
\begin{eqnarray*}\label{LinQ4}
\left(\begin{array}{c} \dot Q\\  \dot p_1 \\ \dot p_2
  \end{array} \right)=M \left(\begin{array}{c}  Q\\  p_1\\ p_2
  \end{array} \right), \quad  M=\begin{pmatrix}
   0 & 1 & 0 \\
  0& 0  & -1\\
c'(x) V& c(x) & 0\\
\end{pmatrix}.
 \end{eqnarray*}
  subject to the initial data
\begin{eqnarray*}\label{cdQ4}
\left(\begin{array}{c}
  Q(0)\\  P(0)
  \end{array} \right)=
  \left(\begin{array}{c}
   1\\  W_0
  \end{array} \right), \quad W_0= \left(\begin{array}{c}  v(x_0) \\  e(x_0)
  \end{array} \right).
  \end{eqnarray*}
  Thus, if in a point $t_*>0$ we have $Q(t_*)=0$, then the solution of \eqref{char_der} blows in the point $t_*$.

  One can check that the system for $(Q, p_1, p_2)$ can be reduced to one linear equation for $Q(t)$:
  \begin{equation}\label{Q}
  \dddot Q + c(x(t)) \dot Q +c'(x(t)) V(x(t)) Q=0,
  \end{equation}
  with the initial conditions
\begin{equation}\label{Q_CD}
  Q(0)=1, \quad Q'(0)=V'_0(x_0), \quad  Q''(0)=-E'_0(x_0),
  \end{equation}
found from \eqref{char_der_CD}.

For arbitrary initial data and an arbitrary doping profile $c(x)$ the system \eqref{char_sol},  \eqref{Q} can be solved numerically with the initial data \eqref{char_sol_CD},  \eqref{Q_CD} for $x_0\in \mathbb R$. This procedure can be performed in reasonable increments with respect to $x$ to cover any desired segment  $[x_-, x_+]\in \mathbb R$.
Note that it is  much easier to determine numerically whether some component of the solution vanishes than to investigate whether this component of the solution goes to infinity or not. We provide examples of such calculations in Sec.\ref{S6}.

\section{Small oscillations}\label{S5}
Now we are going to study how the characteristics of the doping profile influence the properties of small oscillations.
 In this section we assume that $c(x)$ is sufficiently smooth (in fact we only need $C^2$ - smoothness).
\subsection{Change of the period}\label{S51}
\begin{proposition}\label{P1}
The period of the characteristic trajectory $x(t)$ of system \eqref{char_sol}, starting from a point $x_0\in\mathbb R$, depends on the properties of the doping profile $c(x)$. Namely, the following asymptotics holds for the deviation $ \|(V_0(x_0), E_0(x_0))\|$ of order $\varepsilon\ll 1$ from
the origin:
\begin{eqnarray}\label{T}
  T(x_0)&=&\frac{2\pi}{\omega} \left(1+ \Omega \varepsilon^2 + o(\varepsilon^2)    \right),\\
   \omega&=&\sqrt{c(x_0)},\quad
  \Omega=\frac{1}{48} \frac{5 (c'(x_0))^2 - 3 c(x_0)c''(x_0) }{(c(x_0))^3}.\nonumber
\end{eqnarray}
\end{proposition}

\proof
Let us choose $V_0(x_0)=\varepsilon$ as a small parameter. System \eqref{char_sol} implies
\begin{equation}\label{V}
\ddot V +c(x(t)) V=0,\quad \dot x(t)=V,
\end{equation}
therefore if we put $V(0)=\varepsilon$, $\dot V(0)=- E(0)=0$, then the solution can be expanded by the Taylor formula in $\varepsilon$ up to the third order.
 We use the Lindstedt-Poincar\'e method of stretching of the
independent coordinate to avoid secular terms in regular asymptotic expansions
(e.g. \cite{15}, Sec.3) and   set
\begin{equation}
t = s (1 + \varepsilon w_1 + \varepsilon^2 w_2 +\varepsilon^3 w_3 + o(\varepsilon^3)).
\end{equation}
We expand the solution as
\begin{eqnarray*}
 V(s) &=& \varepsilon V_1(s) + \varepsilon^2 V_2(s) +\varepsilon^3 V_3(s) + o(\varepsilon^3),\\
 x(s) &=& x_0+\varepsilon x_1(s) + \varepsilon^2 x_2(s) +\varepsilon^3 x_3(s) + o(\varepsilon^3)\\
 c(x(s))&=&c(x_0)+ \sum\limits_{j=1}^2 \frac{1}{j!}\,c^{j}(x_0) (x(s)-x_0)^j.
  \end{eqnarray*}
From \eqref{V} we get
\begin{eqnarray}\label{exp1}
V_1(s)& =&\cos \omega s , \quad x_1(s)=\frac{\sin \omega s}{\omega}, \\
V_2(s) &=& \frac{-2\sin \omega s c'(x_0) +c'(x_0) \sin 2\omega s-6 \omega^4 s w_1 \sin \omega s}{ 6 \omega^3},\nonumber\\
 x_2(s)& =&\frac{(4 c'(x_0) \cos\omega s -c'(x_0) \sin 2\omega s-3 c'(x_0)+12 \omega^4 s w_1 \cos\omega s)}{12 \omega^3},
\label{exp2}
\end{eqnarray}
therefore $w_1=0$, and
\begin{eqnarray}\label{exp3}
V_3(s)& =& a\omega s \sin \omega s + A_1 \cos \omega s + A_2 \cos 2 \omega s + A_3 \cos 3 \omega s,\nonumber\\
x_3(s)& =& -a s \cos \omega s + B_1 \sin \omega s + B_2 \sin 2\omega s + B_3 \sin 3\omega s,
\end{eqnarray}
with coefficients that depends on $w_2$, $\omega$, $c'(x_0)$, $c''(x_0)$. Computations show that if $a=0$, then $w_2=\Omega$.
$\Box$

\begin{remark}
It is easy to check that $\Omega=0$ for the profile $c(x)=(C_1 x + C_2)^{-\frac32}$, $C_1,\, C_2$ are constant. However, this function is not continuous for all $x\in \mathbb R$. Even if we limit our consideration to the semiaxis, taking into account the following terms of the expansion we will obtain correctors for $T(x_0)$ of order higher than $\varepsilon^2$.
\end{remark}

\subsection{Small perturbations of the zero equilibrium}\label{S52}

We are going to prove the following theorem, which, although a corollary of Theorem \ref{T1}, provides an independent way to show that any small perturbations of the zero equilibrium  necessarily blow up for any non-constant doping profile. 

\begin{theorem}
An arbitrary small smooth perturbation of the zero solution $V=E=0$ of  system \eqref{2} blows up in a finite time for any non-constant smooth doping profile $c(x)>0$.
\end{theorem}

In other words, if there is a point $x_0\in \mathbb R$ such that $c'(x_0)\ne 0$, then a solution with data $(V_0(x_0), E_0(x_0))$ such that $ \|(V_0(x_0), E_0(x_0))\|\ll 1$, blows up in a finite time along the characteristic, starting from $x_0$.
Let us choose $V_0(x_0)=\varepsilon$ as a small parameter and use the results of Proposition \ref{P1}.

Then we can substitute \eqref{exp1}, \eqref{exp2}, \eqref{exp3} to \eqref{Q} as coefficients and apply Floquet's theory (e.g. \cite{Chicone}, section 2.4), which has been successfully used to study other problems, associated with the repulsive Euler-Poisson system. (e.g.  \cite{RD}, \cite{RT}).
The Floquet theory deals with linear systems and equations with periodic coefficients.
According to this theory, for the fundamental matrix $\Psi (t)$ ($\Psi (0)={\mathbb I}$, where ${\mathbb I}$ is the identity matrix) there is a constant matrix $ M$, possibly with complex coefficients such that $\Psi (T)=e^{T M}$, where $T$ is the period of the coefficients. The eigenvalues of the matrix of monodromy $e^{T M}$ are called the characteristic multipliers of the system. If among the characteristic multipliers there are those whose modulus is greater than one, then the zero solution of the linear system under study is Lyapunov unstable (\cite{Chicone}, Theorem 2.53). 

In the case of equation \eqref{Q}, the matrix of monodromy $e^{T M}$ has dimension $(3\times 3)$ and calculating the eigenvalues is quite a difficult task even with the help of a computer algebra package. To simplify our work we integrate \eqref{Q} with respect to $t$ and get (taking into account \eqref{Q_CD})
\begin{equation}\label{Q2}
  \ddot Q(t)+c(x(t)) Q(t)=c(x_0)-E'(x_0)=\rm const,
\end{equation}
\begin{equation}\label{Q2CD}
  Q(0)=1, \quad Q'(0)=V'_0(x_0).
  \end{equation}
Let us study the homogeneous part of \eqref{Q2} (it coincides with \eqref{V}, so we can use the results of Proposition \ref{P1}). In particular, the values for $x_k(t)$, $k=1, 2, 3$ can be taken from \eqref{exp1}, \eqref{exp2}, \eqref{exp3} with the change of $s$ to $t$ and with $w_k=0$.  

The fundamental matrix can be found in the form of the expansion $\Psi(t)=\Psi_0(t)+ \Psi_1(t) \varepsilon + \Psi_2(t) \varepsilon^2+ o(\varepsilon^2),$ $ \Psi_0 (0)=\mathbb I$, $\Psi_k(0)=0$, $k\in \mathbb N$. In our case
\begin{eqnarray*}
\Psi(t)= \begin{pmatrix}
   Q_1(t) & Q_2(t)\\
 \dot Q_1(t) & \dot Q_2(t)
\end{pmatrix}
 \end{eqnarray*}
where $Q_1(t)$, $Q_2(t)$ is the fundamental system of solutions, such that
\begin{eqnarray*}
Q_1(0)=1, \, \dot Q_1(0)=0,\\
Q_2(0)=0, \,\dot Q_2(0)=1.
\end{eqnarray*}
We are looking for $Q_i(t)$, $i=1,2$ as an expansion in $\varepsilon$. At each step, we have to solve a linear inhomogeneous system with constant coefficients; the period $T(x_0)$ (see \eqref{T}) depends on the starting point of the characteristic.

Let us denote  $\Psi_k(t)$ the fundamental matrix constructed for the truncated expansion $\bar Q_{ik}(t)=Q_{i0}(t)+ Q_{i1}(t) \varepsilon +\cdots +Q_{ik}(t) \varepsilon^k$, $i=1,2$. Further, we denote $T_k$ the period of $c(x(t))$, expanded by the Taylor formula in the neighbourhood of $x_0$  with the substitution of $x_1(t), \dots, x_k(t)$.  The eigenvalues of the matrix $\Psi_k(T_k)$ are denoted as  $\lambda_{ik},$ $i=1,2$.  Calculating the eigenvalues at each step  can be done analytically using a computer algebra package (for example, MAPLE).
To prove instability, we need to find an expansion up to  order  $k$ in $\varepsilon$ such that among $ \lambda_{ik}$ there is  a value greater than one in absolute value.

It can be readily shown  that $ \lambda_{i 0}= \lambda_{i 1}=1$.  However,  for $k=2$ we get
\begin{eqnarray*}
\lambda_{\pm 2}=1 \pm (c'(x_0))^{1/2} \left(\frac{2\pi \Omega} {\omega}\right)^{3/2} \, \varepsilon^{7/2},\\
\end{eqnarray*}
see \eqref{T} for the notation. Here we write $\lambda_{\pm 2}$ instead of $\lambda_{i 2}$, $i=1,2$.
  We  see that for any combination of signs of $\varepsilon$, $\Omega$ and $c'(x_0)$ there always there exists an eigenvalue greater than one in absolute value.
  
  If the fundamental system of the homogeneous part \eqref{Q2} contains an oscillating function with an exponentially increasing amplitude, then the solution to the Cauchy problem \eqref{Q2}, \eqref{Q2CD} has the same behavior and there is a point $ t_*>0$ such that $ Q(t_*)=0$.
  The theorem is proved. $\Box$

\begin{remark} We can see that the value
\begin{equation}\label{inst}
m(x)=\Big|c'(x) \left(\frac{ \Omega(x)} {\omega(x)}\right)^{3}\Big|
\end{equation}
can be considered as a measure of instability for the specific doping profile for the small perturbations of the equilibrium state. Thus, to find the most "dangerous" point, we have to find the point $x=m_*$, the maximum of $m(x)$.
\end{remark}

\section{Numerical study}\label{S6}

For the numerical examples we choose the initial data in the form of a standard laser pulse \cite{esarey09}, \cite{CH18}:
\begin{equation}\label{DE}
E_0= a \, x \, e^{-\frac{x^2}{2}}, \qquad V_0=0.
\end{equation}

\medskip

1. The first test is the case of the constant background profile $c(x)=1$.
From \eqref{Q}, \eqref{Q_CD} we find
$$Q(t)=1-2E'_0(x_0)\sin^2 \frac{t}{2},$$
and for $E'_0(x_0)<\frac{1}{2}$ the function  $Q(t)$ does not vanish at all, and therefore the solution does not blow up along the characteristic starting from $x_0$, what corresponds to  criterion \eqref{crit2}. For the initial data \eqref{DE} this condition holds for all $x_0\in \mathbb R$ for $a<\frac12$. Otherwise, if $E'_0(x_0)\ge\frac{1}{2}$ for a point $x_0$, then the solution blows up within a time $t_*<2\pi$.

\medskip

2. The second example is
\begin{equation}\label{Ex2}
c(x)=1+  \frac{0.3}{1+x^2}.
\end{equation}
In the data \eqref{DE} we choose  $a=0.3$, so for any {\it constant} $c(x)=C \in [1, 1.3]$ the solution of the Cauchy problem \eqref{2}, \eqref{DE} does not blow up.

We present the dependence of the blow-up time on $x_0$ for $x_0\in [0.3, 1.3]$ computed from \eqref{char_sol},  \eqref{Q}, \eqref{char_sol_CD},  \eqref{Q_CD} by the Fehlberg fourth-fifth order Runge-Kutta method with the space step $h=0.01$. The function $Q(t)$ oscillates with an increasing amplitude,
the blow-up time is about $75$, it attains on the characteristic starting from the points $x_0 \approx \pm 0.6$. See Fig.1, left.  For the starting points $x_0\in (0, 0.3)$ the blow-up time is greater than $130$. The picture is symmetrical about the origin.
For this profile, the maxima of  $m(x)$ found from \eqref{inst} are approximately $0.2$ and $1$, it is different from the minimum  point of blow-up.
The initial profile \eqref{DE} cannot be considered as a small perturbation of the zero equilibrium, therefore we can conclude that for this case the shape of the initial data does not have a significant effect on the position of the initial point $x_0$ of the characteristic along which the derivatives go to infinity first, which determines the blow-up time of the entire solution of the Cauchy problem (i.e. the minimum of the blow-up time over all $x_0\in \mathbb R$).

 The stepwise nature of the dependence of the blow-up time on the initial point for a specific characteristic is explained by the oscillating nature of $Q(t)$. Namely, a jump in the graph occurs when $Q(t)$ crosses the $Q=0$ axis during the next oscillation.

3. The third example is for the periodic symmetric profile
\begin{equation}\label{Ex3}
c(x)=1+ 0.3 \cos 2 x.
\end{equation}
In the data \eqref{DE} we choose again $a=0.3$ and present the computations on the period $T=\pi$. 
  Near the points, where  $c'(x)=0$ (i.e. $x=0,\, \frac{\pi}{2},\, \pi$) the derivatives go to infinity within a time  greater than $250$. The blow-up time for the solution to the Cauchy problem is less than $25$. See Fig.1, right.
For this profile the maximal points $m_*$ on the period, found from \eqref{inst} are approximately $1.4$ and $1.7$.

\bigskip

  \begin{figure}[htb]
\begin{minipage}{0.4\columnwidth}
\includegraphics[scale=0.3]{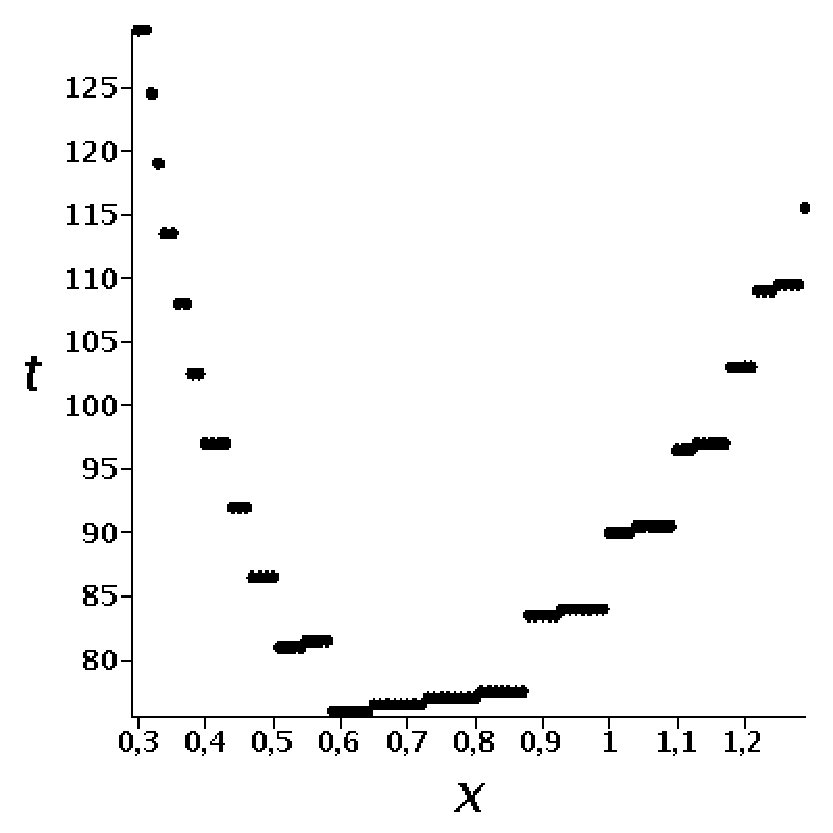}
\end{minipage}
\hspace{1.5cm}
\begin{minipage}{0.4\columnwidth}
\includegraphics[scale=0.3]{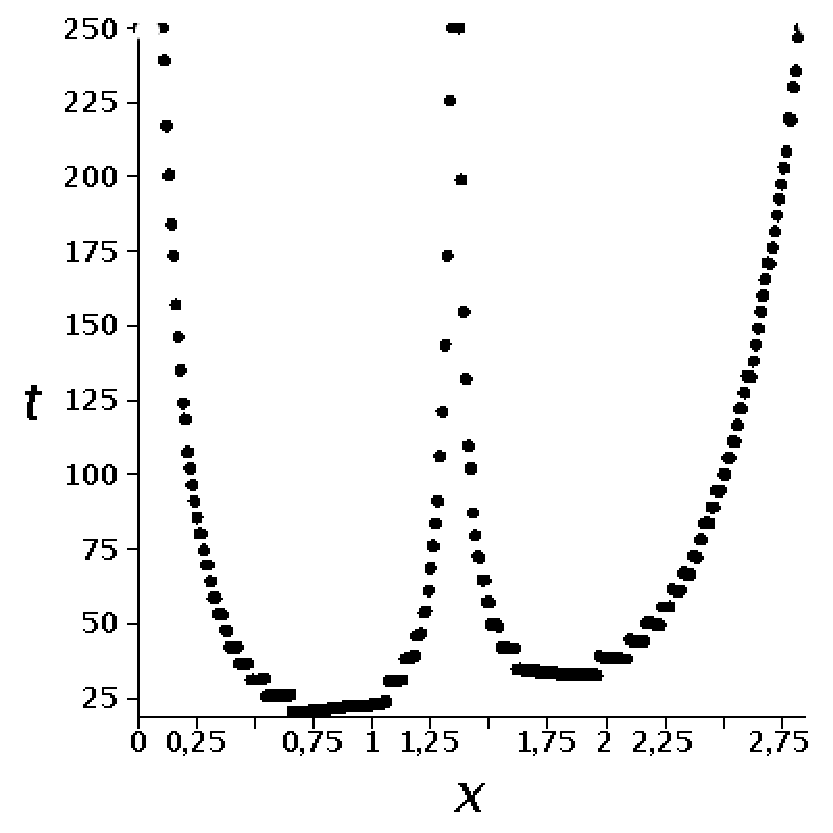}
\end{minipage}

\caption{The blow-up time along characteristics starting from points $x_0$ of the axis $x$ for the initial data \eqref{DE}, $a=0.3$. Left: profile \eqref{Ex2}, right: profile \eqref{Ex3}.}\label{Pic4}
\end{figure}

4. The fourth example is
\begin{equation}\label{Ex4}
c(x)=\frac{1}{(K|x|+1)^2},\quad K>0,
\end{equation}
which is interesting since $m(x)=\rm const$, $x\ne 0$. After example 3, one may get the impression that when determining the time of formation of a singularity, it is not so much the form of the initial data that is important as the doping profile. With this example we show that this is, generally speaking, not the case.
 Fig.2, made for $K=1$ for the interval $x_0\in (0.2, 2.85)$, illustrates that the blow-up time for the solution in this case is minimal for the points closer to the maximum of $s(x)=E_0'(x)$. In the data \eqref{DE} we choose again  $a=0.3$.

\bigskip

  \begin{figure}[htb]
\begin{minipage}{0.4\columnwidth}
\includegraphics[scale=0.3]{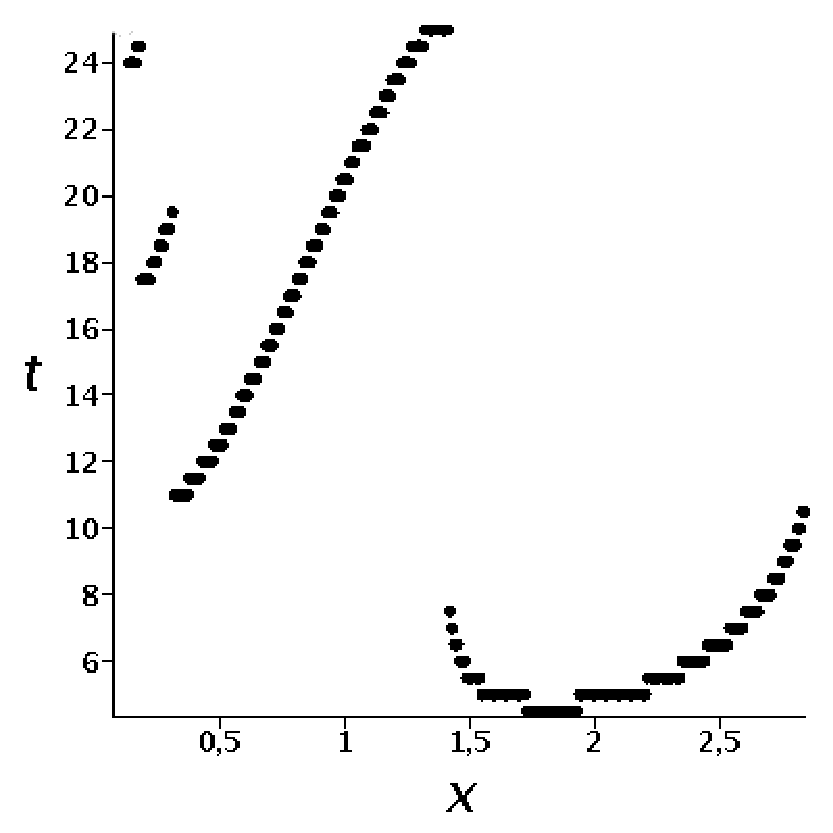}
\end{minipage}
\hspace{1.5cm}
\begin{minipage}{0.4\columnwidth}
\includegraphics[scale=0.3]{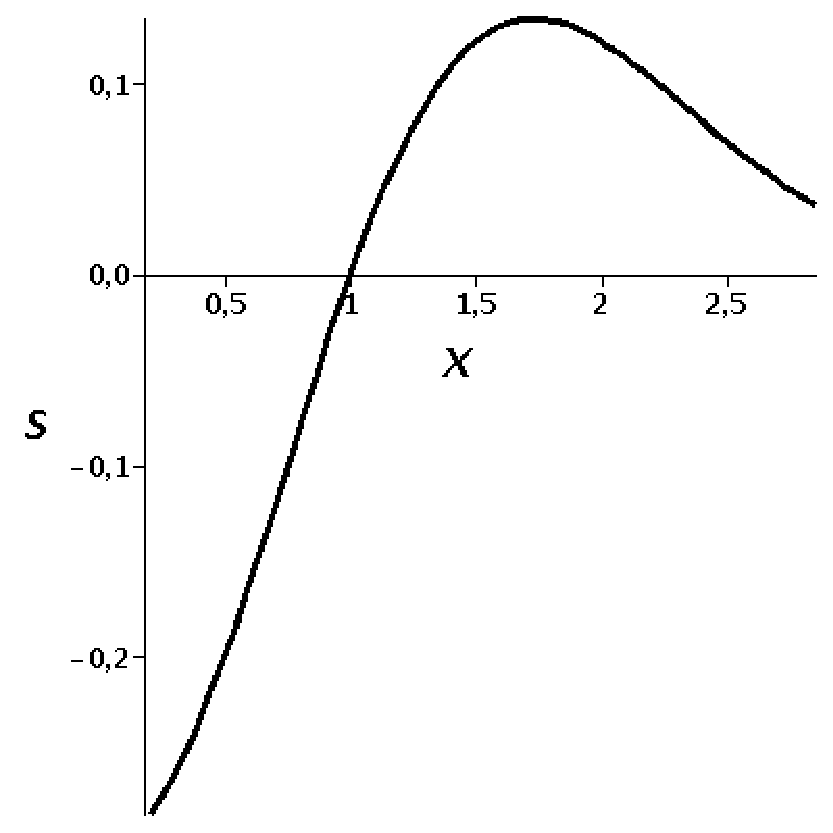}
\end{minipage}

\caption{Left: The blow-up time along characteristics starting from points $x_0$ of the axis $x$ for the initial data \eqref{DE}, $a=0.3$, for the profile \eqref{Ex4}, $K=1$.  Right: the behavior of $s(x)=E_0'(x)$ in the same interval. }\label{Pic4}
\end{figure}

\section{Discussion}\label{S7}

1. Our results, in particular, imply that there is no smooth stationary solution of $(V(x), E(x))$ other than the zero equilibrium. Let us show this independently.
Indeed, from the second equation of \eqref{2} it follows that if $V(x)\ne 0$, then $E_x=c(x)>0$, $\Psi_{xx}= c(x)>0$ and $\Psi$ is strictly concave. Therefore, the identity $\frac12 V^2=-\Psi +\rm const$, following from the first equation of \eqref{2}, cannot be satisfied for all $x\in\mathbb R$. Moreover, this solution has no physical meaning, since for it $n=0$.

Note that the traveling wave solution, depending on the self-similar variable $x-Wt$, $W=\rm const$, does not exist for a variable doping profile, in contrast to the case $c(x) =C=\rm const$. However, even for a constant density background, a smooth traveling wave exists only for sufficiently large $W$ (see \cite{RChZAMP21}).

Note also that for  the attractive Euler-Poisson equations there exists globally smooth solutions for a variable background density $c(x)$, see \cite{Bhat20}, Theorem 2.2, where a periodic $c(x)$ was considered.

\medskip

2. Similar to the model \eqref{1}, we can consider the Euler-Poisson equations with damping
\begin{equation}\label{1A}
\dfrac{\partial V}{\partial t}+ (V \cdot \nabla ) V =-\nabla \Psi - q V,\quad \dfrac{\partial n}{\partial t}+ {\rm div}\,( n V)
=0,\quad \Delta \Psi =c(x) -n.
\end{equation}
In this case, the technique presented in this article allows us to reduce the question on the blow up of a  solution of \eqref{1A} in a finite time to the analysis of the possibility of vanishing of $Q(t)$  in a point $t_*>0$. The function $Q$ solves the equation
  \begin{equation*}\label{QA}
  Q'''+ q Q''+ c(x(t)) Q'+c'(x(t)) V(x(t)) Q=0,
  \end{equation*}
with the initial data \eqref{Q_CD},
$V(t)$ and $x(t)$ are the solutions of
\begin{equation*}\label{VA}
\ddot V + q \dot V + c(x(t)) V=0,\quad \dot x(t)=V,
\end{equation*}
with the initial data
\begin{equation*}\label{char_sol_CDA}
x(0)=x_0, \quad V(0)=V_0(x_0), \quad V'(0)=-E_0(x_0)-q V_0(x_0).
\end{equation*}
Standard calculations show that if the initial data $V_0(x), E_0(x)$ are such that $V_0(x_0), E_0(x_0), V'_0(x_0), E'_0(x_0)$ are sufficiently small in absolute value for each $x_0\in \mathbb R$, then $Q(t)>0$ along each characteristic, and therefore  the respective solution  of \eqref{1A}, with the initial data
\begin{equation*}\label{CDA}
(V,\Psi')|_{t=0}=(V_0(x), \Psi'_0(x))\in {C^2} ({\mathbb R}),
\end{equation*}
keeps the initial smoothness. In other words, in the case $q>0$ there is a small neighborhood of the zero equilibrium in  the $C^1$ - norm such that the solution with initial data from this neighborhood remains smooth. Moreover, as was shown in \cite{RD2024}, this solution asymptotically tends to zero equilibrium as $t\to \infty$.

\medskip

3. For the repulsive Euler-Poisson equations in many spatial dimensions we can consider solutions with radial symmetry and  a variable radially symmetric density background. This problem was considered for the constant density background in \cite{R22_Rad}, where it was found that for $d\ne 1$ and $d\ne 4$ any non-trivial solution other than a simple wave (i.e. $E=E(V)$) blows up. Since simple waves do not exist for a variable doping profile, it is natural to expect that any non-trivial solution blows up.

Note that a thorough analysis of the blow-up conditions was carried out for the radially symmetric case for the attractive case and the case of zero background, as well as for the repulsive case with a constant non-zero background at $d=4$ was made in  \cite{Bhat23}.
 In fact, for the repulsive case the existence of a global smooth solution requires the period of oscillations along every characteristic to be identical. This means that the oscillations are isochronous (i.e. the period does not depend on the amplitude). As follows from \cite{R22_Rad},
 for the case of a constant background, the oscillations are determined by the function $F(t)$, which is a solution of
 a nonlinear Li\'enard type equation
 $$F''+(2+d)\,F\, F' +F+d \, F^3=0,$$
 and the Sabatini criterion \eqref{tau} from Theorem \ref{S} \cite{Sabatini}  implies that the oscillations are  isochronous if and only if $d=1$ and $d=4$. As follows from Lemma \ref{Lem1}, any general nontrivial radially symmetric oscillations blow up for $d\ne 1$ and $d\ne 4$. In \cite{R22_Rad} this is proven analytically for  small perturbations of the zero steady state and numerically for any perturbations. Now we can prove the latter fact analytically.

In \cite{Chen23}, a global weak solution was constructed for the radially symmetric case of the repulsive Euler-Poisson equations with pressure with a variable doping profile.

\section*{Acknowledgements}

 Supported by RSF grant 23-11-00056 through RUDN University.


\begin{thebibliography}{99}


\bibitem{Bhat20}    M. Bhatnagar, H. Liu,
Critical thresholds in 1D pressureless Euler-Poisson systems with variable background, Physica D: Nonlinear Phenomena,
{\bf 414}, 132728 (2020).


\bibitem{Bhat23}    M. Bhatnagar, H. Liu,
A complete characterization of sharp thresholds to
spherically symmetric multidimensional pressureless
Euler-Poisson systems, arXiv:2302.04428 (2023).

\bibitem{Carillo}J.A. Carrillo, R. Shu. Existence of radial global smooth solutions to the pressureless Euler-Poisson equations with quadratic confinement. Arch Rational Mech Anal {\bf 247}, 73 (2023).

\bibitem{Chen23}
G.-Q. G. Chen, L.He, Y.Wang, D.Yuan,
Global solutions of the compressible Euler-Poisson equations for plasma with doping profile for large initial data of spherical symmetry,
arXiv:2309.03158 (2023).

\bibitem{Chicone} Chicone C., {\em Ordinary Differential Equations with Applications}, Springer-Verlag: New York, 1999.

\bibitem{CH18} Chizhonkov E.V., {\em Mathematical aspects of modelling oscillations and wake waves in plasma}, CRC Press, 2019.




\bibitem{RD} M. I. Delova, O. S. Rozanova,
 {\it The interplay of regularizing factors in the model of upper hybrid oscillations of cold plasma}, Journal of Mathematical Analysis and Applications, \textbf{515}: 2 (2022).



\bibitem{Riccati} G. Freiling,
A survey of nonsymmetric Riccati equations, Linear Algebra and its
Applications 351-352, 243-270 (2002).


\bibitem{GR75} V. L. Ginzburg, { Propagation of electromagnetic waves in plasma, Pergamon, New York, 1970}.


\bibitem{Guo17}  Y. Guo, L. Han, J. Zhang. Absence of shocks for one dimensional Euler-Poisson system. Arch.
Rational Mech. Anal., {\bf 223}, 1057-1121 (2017).

\bibitem{ELT}S.Engelberg,
H.Liu, E.Tadmor, Critical Thresholds in Euler-Poisson Equations,
Indiana University Mathematics Journal, {\bf 50}, 109-157 (2001).



\bibitem{esarey09}E. Esarey, C. B. Schroeder, and W. P. Leemans, { Physics of laser-driven plasma-based
electron accelerators,  Rev. Mod. Phys.},  {\bf 81}(2009),
1229-1285.

\bibitem{15}
A. H. Nayfeh, Perturbation methods, John Wiley \& Sons, New York, 2000.

\bibitem{Radon} {W. T. Reid,} Riccati Differential Equations, Academic Press, New York, 1972.


\bibitem{RChZAMP21}
O.S. Rozanova, E.V. Chizhonkov, { On the conditions for the breaking
of oscillations in a cold plasma, { Z. Angew. Math. Phys.},
 }{\bf 72} (2021), 13. 






\bibitem{R22_Rad}
O.S. Rozanova, {\it On the behavior of multidimensional radially symmetric solutions of the repulsive Euler-Poisson equations}, Physica D: Nonlinear Phenomena  {\bf 443}, 133578 (2023). 

\bibitem{RT}  O.S. Rozanova, M.K. Turzinsky, On the properties of affine solutions of cold plasma equations, Communications in Mathematical Sciences {\bf 22},  215-226 (2024).

\bibitem{RD2024} O.S. Rozanova, M.I. Delova, On radially symmetric oscillations of a collisional cold plasma, Mathematical Methods in the Applied Sciences, {\bf 47}(11),  8385-8399 (2024).





\bibitem{RY}B. L. Rozhdestvensky, N. N.Yanenko,  Systems of Quasi-Linear Equations. Am. Math. Soc. Monograph 55 (1983).

\bibitem{Sabatini}
M. Sabatini,  On the period function of Li\'enard systems. J. Differ. Equ. {\bf 152}, 467-487 (1999).

\bibitem{Yuen}
M. Yuen. Blowup for the Euler and Euler-Poisson equations with repulsive forces. Nonlinear Analysis
{\bf 74}, 1465-1470 (2011).

\end{thebibliography}
\end{document}